\documentclass[journal]{IEEEtran}
\ifCLASSINFOpdf
\else
\fi
\usepackage{graphicx} 
\usepackage{amsmath} 
\usepackage{amssymb}  
\usepackage{amsthm}
\usepackage{mathtools}
\usepackage{xcolor}
\usepackage{cite}
\usepackage{caption}
\newtheorem{theorem}{Theorem}

\newtheorem{lemma}[theorem]{Lemma}
\newtheorem{corollary}[theorem]{Corollary}
\newtheorem{definition}{Definition}{}
\usepackage[ruled,noend]{algorithm2e}

\usepackage[colorlinks=true,allcolors=steelblue]{hyperref}
\definecolor{steelblue}{RGB}{70,130,180}
\title{\LARGE \bf
Generalized Outer Bounds on the Finite Geometric Sum of Ellipsoids
}

\author{Navid Hashemi, Justin Ruths
 \thanks{}
\thanks{The authors are with the Department of Mechanical Engineering at the University of Texas at Dallas, 800 W. Campbell Rd, Richardson, TX, USA. Email: {\tt\small (navid.hashemi, jruths)@utdallas.edu}}%
}
 
\begin{document}

\maketitle
\thispagestyle{empty}
\pagestyle{empty}

\begin{abstract}
General results on convex bodies are reviewed and used to derive an exact closed-form parametric formula for the boundary of the geometric (Minkowski) sum of $k$ ellipsoids in $n$-dimensional Euclidean space. Previously this was done through iterative algorithms in which each new ellipsoid was added to an ellipsoid approximation of the sum of the previous ellipsoids. Here we provide one shot formulas to add $k$ ellipsoids directly with no intermediate approximations required. This allows us to observe a new degree of freedom in the family of ellipsoidal bounds on the geometric sum. We demonstrate an application of these tools to compute the reachable set of a discrete-time dynamical system.
\end{abstract}

\section{INTRODUCTION} \label{sec_intro}

Ellipsoidal bounds are efficient representations of sets as they have various appealing properties and require the tuning of relatively few parameters compared with other representations. The convenience of ellipsoids, however, also brings conservatism through the over-approximation of exact sets. For this reason methods tend to seek to identify the most representative ellipsoidal bounds, by finding optimal ellipsoidal bounds by minimizing, for example, the volume. This task is particularly challenged when the exact set is the geometric, or Minkowski, sum of multiple ellipsoids. Although the geometric sum of convex sets is itself a convex set, the geometric sum of several ellipsoids is not, in general, an ellipsoid \cite{lutwak1993brunn,lutwak1996brunn}. Motivated by the ellipsoid sums that arise in the calculation of reachable sets in dynamical systems, in this paper we provide a parameterized expression for the exact boundary of the geometric sum of $k$ ellipsoids as well as a closed form solution for an outer ellipsoidal bound.


Although it is straightforward to define the concept of the geometric sum, it is challenging to quantify the result in closed form to be used, for example, as part of a design problem. There are a variety of results in the literature on the characterization of the geometric sum and difference for two or three dimensional ellipsoids, e.g., see \cite{behar2011dynamic}. In addition the geometric sum has been deeply studied for polygons (two dimension) and polyhedra (three dimension), e.g., see \cite{guibas1983kinetic,agarwal2002polygon,fogel2007exact,hachenberger2009exact,hartquist1999computing}. 
The related algorithms are mainly based on either the computation of the convolution of geometric boundaries \cite{guibas1983kinetic} or polygon/polyhedra decompositions \cite{agarwal2002polygon,fogel2007exact,hachenberger2009exact,hartquist1999computing}. Sums of curved regions/surfaces have also been studied \cite{bajaj1989generation,kaul1995computing,lee1998polynomial,sack1999handbook}. Here in this work, we provide novel one shot proofs to recover known results in the systems and control literature \cite{kurzhanskiui1997ellipsoidal} regarding outer bounding ellipsoids of geometric sums of ellipsoids. The new perspective offered by these proofs allows us to generalize several of these existing results. We specifically show the existence of a new term (degree of freedom) in parameterizing the family of ellipsoids that provide tight outer bounds on the actual geometric sum.

We leverage the proposed results to compute outer ellipsoidal approximations for the reachable set of a discrete-time linear control system subject to ellipsoidally-bounded input sets. Such calculations, to date, have been done iteratively and require calculating an outer ellipsoid approximation after the addition of each ellipsoid \cite{kurzhanskiui1997ellipsoidal}. Here, we propose an approach which computes the geometric sum all in one step, which yields a closed-form characterization of the ellipsoid bound and also reveals the opportunity to generalize the existing results.

\section{Geometric Sum}

In this work, we demonstrate how to specify the exact boundary of the geometric sum of multiple ellipsoidally bounded sets. 

\begin{definition}{\cite{ellipsoidal_toolbox}}
The geometric (or Minkowski) sum of two convex sets $\mathcal{S}_1,\mathcal{S}_2\subset\mathbb{R}^n$ is
\begin{equation}
    \mathcal{S}_1\oplus\mathcal{S}_2 = \left\{s_1+s_2\ \Big|\ s_1\in\mathcal{S}_1,\, s_2\in\mathcal{S}_2 \right\}.
\end{equation}
\end{definition}
Note that while the geometric sum of two ellipsoids is a convex set, it is in general not an ellipsoid.

\subsection{Exact Boundary of a Geometric Sum}\label{sec:geometricsum}

\begin{figure}[t]
\centering
\includegraphics[width=0.8\linewidth]{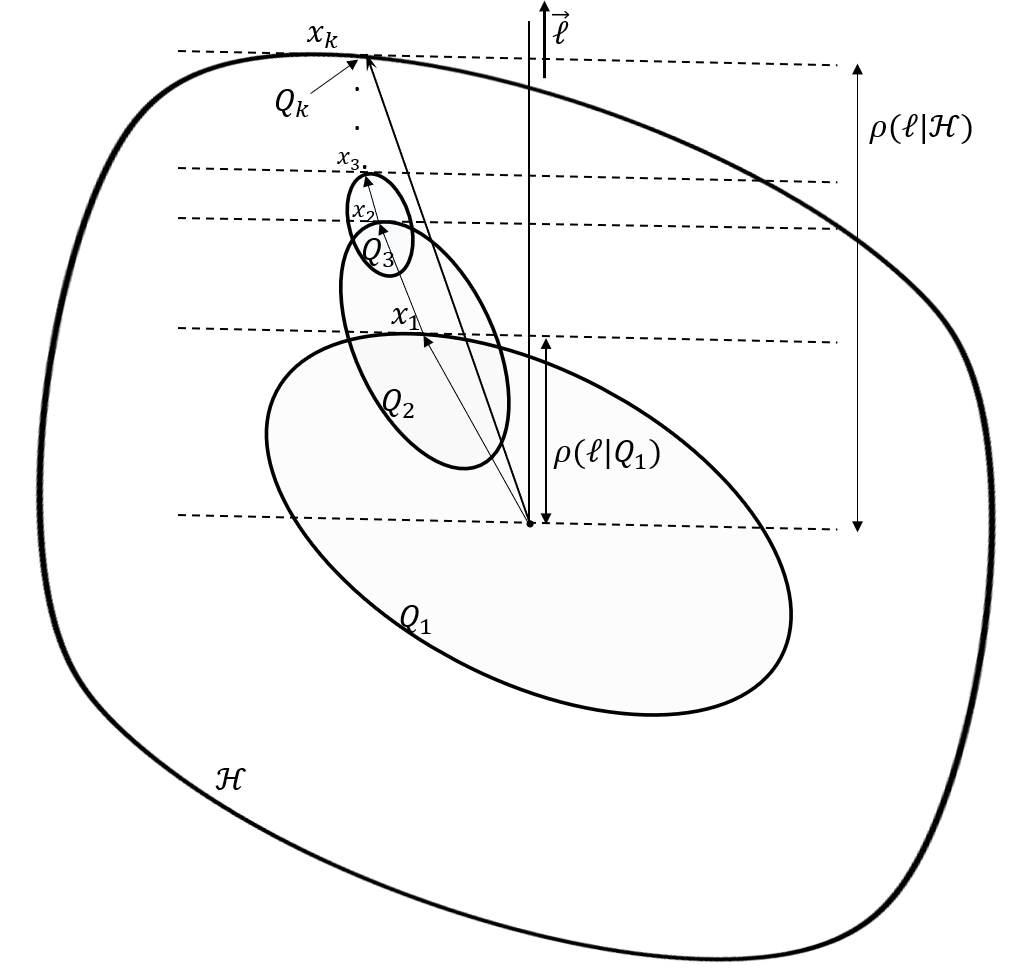}
\caption{Convex set $\mathcal{H}$ as a geometric sum of $k$ ellipsoids. Here the ellipsoids are decreasing in size as would be the case for a stable dynamical system and ellipsoids generated by \eqref{eq:shapematrixstate2}.} \label{geomsum}
\end{figure}

A central tool in working with geometric operations on convex sets is the support function,
\begin{equation}\label{eq:support}
\rho(\ell|\mathcal{H})=\sup_{x\in \mathcal{H}}\langle \ell,x \rangle,
\end{equation}
which can be interpreted as the largest projection of elements in the convex set $\mathcal{H}$ onto the direction given by the unit vector $\ell$. Containment of one convex set within another convex set, $\mathcal{H} \subset \mathcal{S}$, is exactly captured by the support function of the contained set being less than or equal to the support function of the containing set over all choices of $\ell$,
\begin{equation}
\label{eq:subsetlaw}
\rho(\ell| \mathcal{H}) \leq \rho(\ell | \mathcal{S}), \qquad \forall \ell \in \mathbb{R}^n.
\end{equation}
In addition, the support function of a geometric sum of two convex sets can be expressed as the sum of their support functions \cite{kurzhanskiui1997ellipsoidal}. The geometric sum of $k$ ellipsoids centered at zero with shape matrices $Q_1,\dots,Q_k$ (see Fig. \ref{geomsum}) is then
\begin{equation} 
\label{eq:supportsum}
\rho( \ell|\mathcal{H})=\sum_{i=1}^k\rho\big(\ell\,|\,\mathcal{E}(Q_i)\big),
\end{equation}
where the ellipsoid with shape matrix $Q$ is given by
\begin{equation} \label{eq:ellipsoid}
\mathcal{E}(Q) = \left\{ \xi\  | \ \xi^T Q^{-1} \xi\leq 1 \right\}.
\end{equation}
When the convex sets of interest are ellipsoids, the support function can further be expressed as
\begin{equation}\label{eq:ellipscase}
\rho\big(\ell\,|\,\mathcal{E}(Q)\big)=\langle \ell,Q\ell \rangle^{\frac{1}{2}}.
\end{equation}
Therefore, the geometric sum \eqref{eq:supportsum} can be written as
\begin{equation} 
\label{supportsumellips}
\rho(\ell|\mathcal{H})=\sum_{i=1}^k\langle \ell,Q_i\ell \rangle^{\frac{1}{2}}.
\end{equation}

The following result provides a closed form formula to calculate the exact boundary of the geometric sum of a finite collection of ellipsoids based on their support functions. Existing approaches approach the sum of a finite collection of ellipsoids in an iterative fashion, in which the $j$-th ellipsoid is geometrically summed with the result of the sum of ellipsoids $1$ through $j-1$. Although this iterative approach yields the same results in some cases, we illustrate in this paper that our proposed ``one-shot'' approach provides additional insight. Throughout the paper we will use the notation $\overline{1,k}=\{1,2,\dots,k\}$.
\begin{theorem} \label{exactbound} Consider the set of $k$ ellipsoids characterized by the positive definite matrices $Q_i \in \mathbb{V}^{n \times n}$, $i\in\overline{1,k}$. The exact convex shape of their geometric sum is, $\mathcal{H} = \bigoplus_{i=1}^k \mathcal{E}(Q_i)$ with boundary
\begin{equation}\label{eq:exactreachableset}
    \partial\mathcal{H} = \left\{x \ \bigg|\  x=\sum_{i=1}^{k}\langle \ell,Q_i \ell \rangle^{-\frac{1}{2}}Q_i \ell,\ \forall \ell \in \mathcal{B}_1(0) \right\},
\end{equation}
where $\mathcal{B}_1(0)$ is the unit ball in $\mathbb{R}^n$.
\end{theorem}

\textit{Proof:}
One way to interpret the support function \eqref{eq:support} is that $\ell$ provides the normal vector of the hyperplane that is tangent to the convex set $\mathcal{H}$ at the point $x$. By the definition of the geometric sum, the boundary point $x$ can be written as the sum of elements $x_i\in\mathcal{E}(Q_i)$ for $i\in\overline{1,k}$,
\begin{equation}\label{eq:approach}
\textstyle x=\sum_{i=1}^{k}x_i,
\end{equation}
where each $x_i$ is a boundary point on each individual ellipsoid, i.e., the maximizer of
\begin{equation} \label{eq:supportQr}
    \rho\big(\ell\,|\,\mathcal{E}(Q_i)\big)=\sup_{\xi\in \mathcal{E}(Q_i)}\langle \ell,\xi \rangle,
\end{equation}
analogous to \eqref{eq:support}. Hence, for each $i\in\overline{1,k}$, $\ell$ provides the (unit) normal vector of the hyperplane that is tangent to the ellipsoid $\mathcal{E}(Q_i)$ at the point $x_i$. The normal vector of the ellipsoid can also be expressed as the gradient of the level set $x_i^TQ_i^{-1}x_i=1$, hence,
\begin{equation}\label{eq:normalvector}
\textstyle \ell=\frac{\nabla (x_i^T Q_i^{-1} x_i -1)}{\| \nabla (x_i^T Q_i^{-1} x_i -1) \|}=\frac{Q_i^{-1}x_i}{\sqrt{x_i^TQ_i^{-2}x_i}}.
\end{equation}
If we use again the fact that $x_i^TQ_i^{-1}x_i=1$,
\begin{align}
\ell&=\textstyle \frac{(x_i^T Q_i^{-1} x_i) Q_i^{-1} x_i }{ \sqrt{ x_i^T Q_i^{-2} x_i}}
= x_i^T \frac{Q_i^{-1} x_i }{ \sqrt{ x_i^T Q_i^{-2} x_i}}  Q_i^{-1} x_i \nonumber\\
&= x_i^T \ell Q_i^{-1} x_i. \label{eq:normalvector2}
\end{align}
Note that as the maximizer of \eqref{eq:supportQr}, the projection of $x_i$ onto $\ell$ gives the support function,
\begin{equation}\label{eq:supportresult}
    x_i^T\ell=\rho(\ell\, |\, \mathcal{E}(Q_i))=\langle \ell,Q_i \ell \rangle^{\frac{1}{2}}.
\end{equation}
Substituting \eqref{eq:supportresult} into \eqref{eq:normalvector2} yields
\begin{equation} \label{eq:normalvector3}
\ell = \langle \ell,Q_i \ell \rangle^{\frac{1}{2}}Q_i^{-1}x_i.    
\end{equation}
Finally, solving for $x_i$
and plugging these into \eqref{eq:approach} completes the proof.
\hspace*{\fill}~$\blacksquare$


\subsection{Ellipsoidal Outer Bounds of the Geometric Sum}
Because the exact boundary specified by Theorem \ref{exactbound} is a function of $\ell$ and because it is in general not an ellipsoid, we are typically interested in computing an outer bound on this exact set, which helps to maintain tractability in many settings. 

Theorem \ref{subset} characterizes the general family of shape matrices that provide an outer ellipsoidal bound on the exact geometric sum, i.e., $\mathcal{H}\subseteq\mathcal{E}(Q)$. This result is similar to Theorem 2.7.1 in \cite{kurzhanskiui1997ellipsoidal}, which computes the geometric sum of $k$ ellipsoids through pairwise sums. The recursive approach in \cite{kurzhanskiui1997ellipsoidal} approximates the geometric sum of two ellipsoids with an ellipsoid at each iteration which, as we will show through the introduction of the matrix $Q_0$, fundamentally restricts the choice of the outer bound. 

\begin{theorem}\label{subset}
Consider the set of $k$ ellipsoids characterized by the positive definite matrices $Q_i \in \mathbb{V}^{n \times n}$, $i\in\overline{1,k}$. The ellipsoid $\mathcal{E}(Q)$ is an outer bound of their geometric sum $\mathcal{H} = \bigoplus_{i=1}^k \mathcal{E}(Q_i)$, i.e., $\mathcal{H}\subseteq\mathcal{E}(Q)$, if
\begin{equation}\label{elipbound}
\begin{aligned}
Q&=Q_0+Q_u,\\
Q_u&=\sum_{i=1}^k Q_i +\sum_{i,j \in \mathcal{I}} p_{ij}Q_i+p_{ij}^{-1}Q_j,
\end{aligned}
\end{equation}
where $p_{ij} > 0\ \forall(i,j)\in\mathcal{I}=\{(i,j)\ |\ 1 \leq i < j \leq k,\  i,j \in \mathbb{N}\}$ and $Q_0$ satisfies the following characteristics:
\begin{enumerate}
\item $Q_0+Q_u>0$,
\item $\langle \ell,(Q_0+Q_u) \ell \rangle^{\frac{1}{2}} \geq \sum_{i=1}^{k} \langle \ell,Q_i \ell \rangle^{\frac{1}{2}},\ \forall \ell\in \mathcal{B}_1(0)$.
\end{enumerate}
\end{theorem} 
\textit{Proof:} To prove that the convex set of the geometric sum is bounded by the ellipsoid, i.e., $\mathcal{H}\subseteq\mathcal{E}(Q)$, it suffices to prove \eqref{eq:subsetlaw} holds between the support functions $\rho(\ell | \mathcal{E}(Q))$ and $\rho(\ell | \mathcal{H})$. We first define a positive coefficient $p_{ij}$ between the shape matrices $Q_i$ and $Q_j$ and in order to avoid redundancy assume $i<j$, which constructs the set $\mathcal{I}$. Given a unit vector $\ell$, we start with the following sum of squares which is always nonnegative
\begin{equation}\label{eq:sumofsquares}
\sum_{i,j \in \mathcal{I}}\left(\sqrt{p_{ij}}\langle \ell,Q_i \ell \rangle^{\frac{1}{2}}-\sqrt{p_{ij}}^{-1}\langle \ell,Q_j\ell \rangle^{\frac{1}{2}}\right)^2 \geq 0.
\end{equation}
Expanding this series results in
$$
\sum_{i,j \in \mathcal{I}} \left(p_{ij}\langle \ell,Q_i\ell \rangle+p_{ij}^{-1}\langle \ell,Q_j\ell \rangle\right) \geq \sum_{i,j \in \mathcal{I}}2\langle \ell,Q_i\ell \rangle^{\frac{1}{2}}\langle \ell,Q_j \ell \rangle^{\frac{1}{2}}.
$$
Now we add $\sum_{i=1}^k \langle \ell,Q_i\ell \rangle$ to both sides and add $\langle \ell,Q_0\ell \rangle$ to the left side of the inequality. Since $Q_0$ satisfies the two conditions presented in the theorem, the inequality still holds,
\begin{align}
&\langle \ell,Q_0\ell \rangle+\sum_{i=1}^k \langle \ell,Q_i\ell \rangle +\sum_{i,j \in \mathcal{I}} p_{ij}\langle \ell,Q_i\ell \rangle+p_{ij}^{-1}\langle \ell,Q_j\ell \rangle \nonumber\\ 
&\geq   \sum_{i=1}^k \langle \ell,Q_i\ell \rangle +\sum_{i,j \in \mathcal{I}}2\langle \ell,Q_i\ell \rangle^{\frac{1}{2}}\langle \ell,Q_j\ell \rangle^{\frac{1}{2}}. \label{eq:sumofsquaresexpanded}
\end{align}
This allows us to complete the square on the right hand side of the inequality. The left hand side, then, becomes the definition of $Q$ in \eqref{elipbound}, such that it is equal to $ \langle \ell,Q\ell \rangle $. Therefore, $ \langle \ell,Q\ell \rangle \geq \big(\sum_{i=1}^k\langle \ell,Q_i\ell \rangle^{\frac{1}{2}}\big)^2$\!, which means,
\begin{equation}\label{eq:subsetinequality}
\underbrace{\langle \ell,Q\ell \rangle^\frac{1}{2}}_{\rho\big(\ell\,|\,\mathcal{E}(Q)\big)} \geq \underbrace{\sum_{i=1}^k\langle \ell,Q_i\ell \rangle^{\frac{1}{2}}}_{\rho(\ell |\mathcal{H})},
\end{equation}
and the underbrace definitions come from \eqref{eq:ellipscase} and \eqref{supportsumellips}. We conclude then that $\mathcal{H} \subseteq \mathcal{E}(Q)$ from \eqref{eq:subsetlaw}.
\hspace*{\fill}~$\blacksquare$

Although it is written in several terms in \eqref{elipbound}, the shape matrix $Q$ is a linear combination of the individual shape matrices $Q_i$, $Q=Q_0+\sum_{i=1}^k\alpha_iQ_i$, with coefficients $\alpha_i$ that are the combination of $p_{ij}$ and $p_{ij}^{-1}$. If the coefficients of the $Q_i$ are too large or if $Q_0$ is positive definite, then it makes $\mathcal{E}(Q)$ a loose outer bound on $\mathcal{H}$. The former occurs when $p_{ij}$ is either very large or very small (so $p_{ij}^{-1}$ is large). Thus to have a tight outer bound, which is tangent to the exact sum $\mathcal{H}$, the challenge is to select the definition of $p_{ij}$ appropriately. This insight is largely already captured in analogous results that exist in the literature \cite{kurzhanskiui1997ellipsoidal}. What is novel here, is the additional degree of freedom provided by $Q_0$. A tempting trivial solution is to select $Q_0=0$, however, we will show in Section \ref{sec:mintrace} that this choice prohibits the minimum trace bounding ellipsoid from being a tight outer bound on the exact geometric sum.
\begin{figure}[t]
\centering
\includegraphics[width=0.8\linewidth]{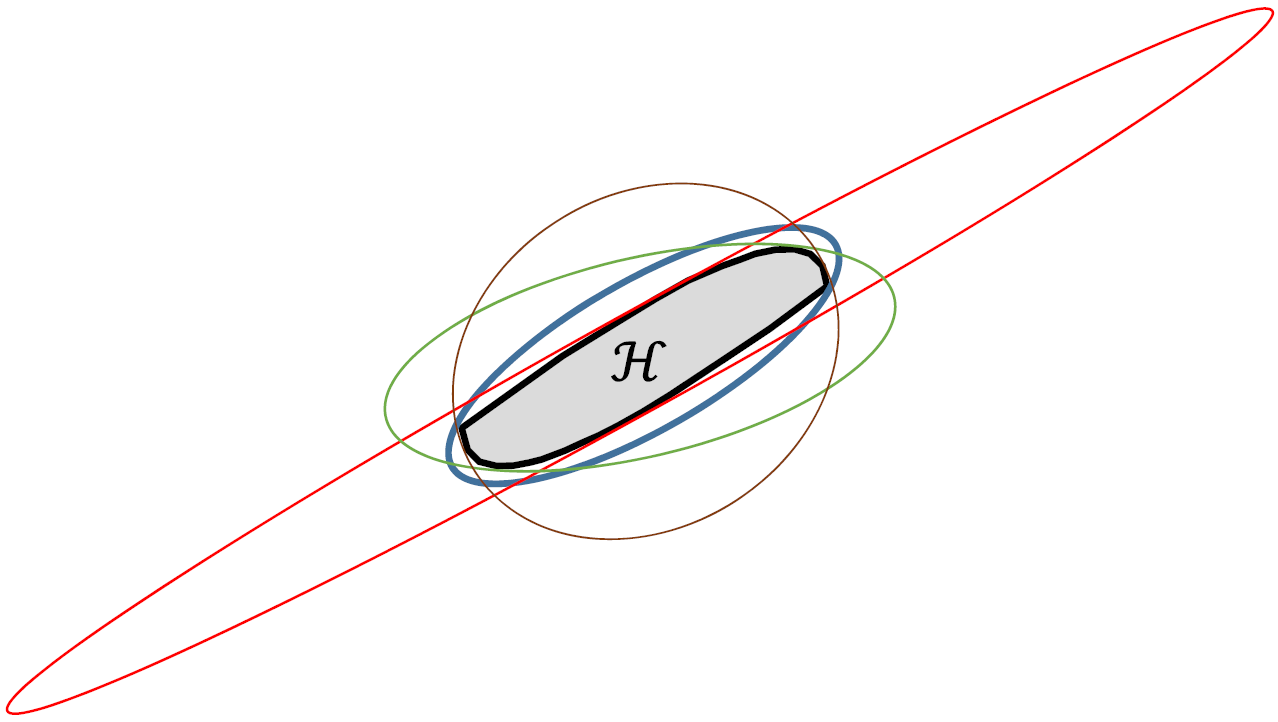}
\caption{There are an infinite number of ellipsoids that are a tight outer bound of a convex set $\mathcal{H}$, each tangent to $\mathcal{H}$ at different points (corresponding to different choices of the unit vector $\ell$). Since many of these tight ellipsoids can actually be quite misleading, we assert that the minimum trace ellipsoid produces the most representative outer bound.
} \label{fig:elipsbound}
\end{figure}

In the current literature, the family of tight (tangent) outer ellipsoidal bounds on the geometric sum is parameterized by a choice of $\ell$ \cite{ellipsoidal_toolbox},
\begin{equation}\label{existant}
    Q=\left( \sum_{i=1}^{k} \langle \ell,Q_i \ell \rangle^{\frac{1}{2}} \right)\left( \sum_{i=1}^{k}\langle \ell,Q_i \ell \rangle^{-\frac{1}{2}}Q_i \right).
\end{equation}
The choice of $\ell$ changes where the outer bound is tangent to the actual geometric sum (see Fig. \ref{fig:elipsbound}). Here, in Corollary \ref{othershapematrix} we show that the one-step approach provided by Theorems \ref{exactbound} and \ref{subset} enable a more general formulation of the family of tangent ellipsoid outer bounds than is provided in \eqref{existant}.   

\begin{corollary}\label{othershapematrix}
Consider the set of $k$ ellipsoids characterized by the positive definite matrices $Q_i \in \mathbb{V}^{n \times n}$, $i\in\overline{1,k}$. Given a unit vector $\ell\in \mathbb{R}^n$ the ellipsoids $\mathcal{E}(Q)$ with
\begin{equation}\label{eq:shapematrix2}
\begin{aligned}
    Q&=Q_0+Q_\ell,\\
    Q_\ell&=\left( \sum_{i=1}^{k} \langle \ell,Q_i \ell \rangle^{\frac{1}{2}} \right)\left( \sum_{i=1}^{k}\langle \ell,Q_i \ell \rangle^{-\frac{1}{2}}Q_i \right),
\end{aligned}
\end{equation}
provide all possible ellipsoidal bounds that touch the boundary of the geometric sum $\mathcal{H} = \bigoplus_{i=1}^k \mathcal{E}(Q_i)$, i.e., $\mathcal{H}\subseteq\mathcal{E}(Q)$, and are tangent at the respective points
\begin{equation}\label{eq:exactreachableset2}
    x=\sum_{i=1}^{k}\langle \ell,Q_i \ell \rangle^{-\frac{1}{2}}Q_i \ell.
\end{equation}
and where $Q_0$ satisfies the following characteristics:
\begin{enumerate}
\item $Q_0+Q_\ell>0$,
\item $\ell \in \ker (Q_0)$,
\item $\langle \tilde\ell,(Q_0+Q_\ell) \tilde\ell \rangle^{\frac{1}{2}} \geq \sum_{i=1}^{k} \langle \tilde\ell,Q_i \tilde\ell \rangle^{\frac{1}{2}},\ \forall \tilde\ell\in \mathcal{B}_1(0)$.
\end{enumerate}
\end{corollary}

\textit{Proof:} 
This corollary can be derived either from Theorem \ref{exactbound} or \ref{subset}. For completeness we provide both as they provide different insight.

\textit{Based on Theorem \ref{subset}}, the outer bound is tangent to the exact geometric sum if \eqref{eq:subsetinequality} holds with equality for some $\ell$, which in turn requires \eqref{eq:sumofsquares} and \eqref{eq:sumofsquaresexpanded} to hold with equality. The former occurs when the $p_{ij}>0$ coefficients are selected as
\begin{equation}\label{eq:tangentchoice}
p_{ij}=\frac{\langle \ell,Q_j\ell \rangle^\frac{1}{2}}{\langle \ell,Q_i\ell \rangle^\frac{1}{2}}.
\end{equation}
Substituting \eqref{eq:tangentchoice} into \eqref{elipbound} results in \eqref{eq:shapematrix2}. For \eqref{eq:sumofsquaresexpanded} to hold with equality, we also require $\langle \ell,Q_0\ell\rangle=0$, which leads to the new condition on $Q_0$, $\ell\in\ker(Q_0)$.

\textit{Based on Theorem \ref{exactbound}}, given a unit vector $\ell$ the corresponding point on the boundary of the geometric sum is given by \eqref{eq:exactreachableset2}.
For the ellipsoid $\mathcal{E}(Q)$ to touch the boundary of the geometric sum at this point we equate the support function of the outerbounding ellipsoid \eqref{eq:ellipscase} with the support function of the geometric sum \eqref{supportsumellips},
\begin{equation}\label{eq:subsetequality}
\langle \ell,Q\ell \rangle^\frac{1}{2} = \sum_{i=1}^k\langle \ell,Q_i\ell \rangle^{\frac{1}{2}}.
\end{equation}
Using the same approach as in Theorem \ref{exactbound} (see \eqref{eq:normalvector3}), we find the relationship between the tight ellipsoidal bound shape matrix, $Q$, the point of tangency $x$, and the direction $\ell$ as,
\begin{equation}\label{eq:ellstar}
 \ell = \underbrace{\langle \ell,Q \ell \rangle^{\frac{1}{2}}}_{\text{\eqref{eq:subsetequality}}} Q^{-1} \underbrace{x}_{\text{\eqref{eq:exactreachableset2}}}.
\end{equation}
Making the substitutions in the underbraces, multiplying both sides by $Q$, and factoring out $\ell$ yields
\begin{align*}
    Q\ell&= \left[\left( \sum_{i=1}^{k} \langle \ell,Q_i \ell \rangle^{\frac{1}{2}} \right)\left( \sum_{i=1}^{k}\langle \ell,Q_i \ell \rangle^{-\frac{1}{2}}Q_i \right) \right] \ell.
\end{align*}
As a positive definite shape matrix, $Q$ is invertible, which leads to the definition of $Q$ in \eqref{eq:shapematrix2}, where $Q_0$ satisfies the proposed conditions. 
\hspace*{\fill}~$\blacksquare$

Like Theorem \ref{subset}, Corollary \ref{othershapematrix} establishes the existence of a set of symmetric matrices $Q_0$ that modify the shape matrices provided in \eqref{existant}, with the trivial choice $Q_0=0$ resulting directly to \eqref{existant}. While in this work we do not specify $Q_0$ any further than the conditions listed in Theorem \ref{subset} or Corollary \ref{othershapematrix}, it clarifies that the existing characterization of shape matrices omit the additional degree of freedom provided by $Q_0$. 



\subsection{Minimum Trace Ellipsoidal Outer Bound}\label{sec:mintrace}


Figure \ref{fig:elipsbound} easily demonstrates that a tight ellipsoidal bound does not necessarily provide a representative outer bound for a convex set. Even ellipsoids that minimize volume can lead to highly eccentric bounding ellipsoids. The minimum trace ellipsoid tends to produce a more representative outer bound.  The trace of the shape matrix is the sum of the eigenvalues and the lengths of the principal axes of the ellipsoid are the singular values. Thus, by minimizing the trace it penalizes long principal axes. We denote this minimum trace ellipsoid by $\mathcal{E}(Q^*)$. A further benefit of this choice is that facilitates an analytic expression for its shape matrix, which is derived in Lemma \ref{mintrace}.

\begin{lemma}\label{mintrace}
Consider the set of $k$ ellipsoids characterized by the positive definite matrices $Q_i \in \mathbb{V}^{n \times n}$, $i\in\overline{1,k}$. The ellipsoid $\mathcal{E}(Q^*)$ is an outer bound of their geometric sum $\mathcal{H} = \bigoplus_{i=1}^k \mathcal{E}(Q_i)$, i.e., $\mathcal{H}\subseteq\mathcal{E}(Q^*)$, where 
\begin{equation}\label{eq:starformula}
\begin{aligned}
    Q^*&=\left( \sum_{i=1}^k \sqrt{\mathrm{tr}(Q_i)}\right)\left(\sum_{i=1}^k \frac{Q_i}{\sqrt{\mathrm{tr}(Q_i)}}\right).
\end{aligned}
\end{equation}
Specifically, of all the shape matrices that correspond to ellipsoid outer bounds of $\mathcal{H}$ using Theorem \ref{subset} with $Q_0=0$, $Q^*$ has the minimum trace.
\end{lemma}
 
\textit{Proof:} Taking $Q_0=0$ and starting from the general condition on $p_{ij}$ for $\mathcal{E}(Q)$ to be an outer bound, \eqref{elipbound} of Theorem \ref{subset}, we take the trace operation of both sides
\begin{equation}\label{eq:1}
\mathrm{tr}(Q)=\sum_{i=1}^k \mathrm{tr}(Q_i) +\sum_{i,j \in \mathcal{I}} p_{ij}\mathrm{tr}(Q_i)+p_{ij}^{-1}\mathrm{tr}(Q_j),
\end{equation}
making $\mathrm{tr}(Q)$ a function of $p_{ij}$, $(i,j)\in \mathcal{I}$. To find the minima, we find the stationary points of $\mathrm{tr}(Q)$ with respect to $p_{ij}$,
\begin{equation}\label{eq:2}
\frac{\partial{\mathrm{tr}(Q)}}{\partial{p_{ij}}}=\mathrm{tr}(Q_i)-\frac{\mathrm{tr}(Q_j)}{p_{ij}^2}=0.
\end{equation}
The Laplacian is strictly positive due to the positive definiteness of $Q_1,\ldots,Q_k$,
\begin{equation}
  \textstyle \nabla^2 \mathrm{tr}(Q)=\textbf{diag} \left[ \ \   \left\{ \frac{2\mathrm{tr}(Q_j)}{p_{ij}^3} \right\}_{(i,j) \in \mathcal{I}} \ \  \right] >0,
\end{equation}
where the pairs $(i,j)$ are ordered first by $i$ and then by $j$. Hence the solution of \eqref{eq:2} is a unique global minimum,
\begin{equation}\label{tightestchoice}
    \textstyle p_{ij}^*=\sqrt{\frac{\mathrm{tr}(Q_j)}{\mathrm{tr}(Q_i)}}.
\end{equation}
Substituting this choice of $p_{ij}$ into \eqref{elipbound} completes the proof. 
\hspace*{\fill}~$\blacksquare$

The result of Lemma \ref{mintrace}, namely \eqref{eq:starformula}, is a well-known result \cite{kurzhanskiui1997ellipsoidal}. What we have shown in Lemma \ref{mintrace} is that this formula explicitly requires $Q_0=0$ in Theorem \ref{subset}. A smaller trace is achievable by introducing a nonzero $Q_0$ term. This has two implications: (a) that \eqref{eq:starformula} does not find the true minimum trace when $Q_0$ is nonzero, and (b) that \eqref{eq:starformula} is not necessarily tangent to the exact geometric sum. If it was tangent, then the trace could not be reduced further.

Another observation, which for space we do not prove explicitly here, is that in the first proof of Corollary \ref{othershapematrix}, the coefficients $p_{ij}$ were shown to follow \eqref{eq:tangentchoice} if the ellipsoidal bound was to be a tight outer bound. In contrast, the Lemma \ref{mintrace} claims that \eqref{tightestchoice} specifies the same coefficients. It can be shown that the $n$-dimensional unit vector (which has $n-1$ degrees of freedom) cannot be selected to reconcile the two equations for $p_{ij}$ over all $(i,j)\in\mathcal{I}$.

\section{Application in Computing Reachable Sets}

Consider the dynamics of a linear time varying (LTV) system with $r$ independent inputs, 
\begin{equation}
s_{k+1}=A_k s_k+B_{1,k}u_{1,k}+B_{2,k}u_{2,k}+ \cdots + B_{r,k}u_{r,k},
\end{equation}
where $s_k \in \mathbb{R}^n$ represents the state and $u_{i,k} \in \mathcal{E}({R}_{i,k})$, $R_{i,k}\in\mathbb{R}^{n\times m_i}$, $i=\overline{1,r}$ are ellipsoidally bounded (independent) inputs. 
The previous sections of this paper provide tools to compute the exact boundary of the reachable states of this dynamical system and a practical outer bound. We use the geometric sum to explore all states $s_k$ that are reachable through the different realizations of $u_{i,j}$ over $j\in\overline{1,k-1}$.

To compute the reachable set, we express $s_k$ as a function of the prior  inputs. For simplicity, we assume the initial state $s_1=0$, however, incorporating an ellipsoidally bounded set of initial conditions would be a simple extension of this framework. 
\begin{equation}\label{eq:state_seq}
     s_k=\sum_{i=1}^{r} \left(\sum_{j=1}^{k-1}D_jB_{i,j}u_{i,j}\right),\qquad
    D_{j}=\prod_{\kappa=j+1}^{k-1}A_{\kappa}. 
\end{equation}
Equation \eqref{eq:state_seq} shows that the reachable set, the union of all possible state values under any combination of admissible inputs, can be expressed as the geometric sum of transformed, independent, and ellipsoidally bounded inputs $u_{i,j}$.
\begin{lemma}{\cite{ellipsoidal_toolbox}} \label{lem:transformation}
Given vector $x\in\mathcal{E}(Q,c)$, with shape matrix $Q \in \mathbb{V}^{n \times n}$ and mean value $c \in \mathbb{R}^n$, then for any linear mapping $T(x)=Mx$ with $A \in \mathbb{R}^{n \times n}$ the vector $\xi=T(x)=Mx$ lies within an ellipsoid with shape matrix $MQM^T$ and mean value $Mc$, i.e., $\xi\in\mathcal{E}(MQM^T,Mc)$.
\end{lemma}
Using Lemma \ref{lem:transformation} to compute the transformed shape matrices of the input ellipsoids, the geometric sum based on \eqref{eq:state_seq} that is the reachable set is given by
\begin{equation}\label{eq:shapematrixstate2}
\mathcal{R}_{k}=\bigoplus_{i=1}^{r}\bigoplus_{j=1}^{k-1} \mathcal{E}\left(\Sigma_{i,j} \right),
\end{equation}
with $\Sigma_{i,j}=D_j B_{i,j} R_{i,j} B_{i,j}^T D_j^T$.

\subsection{Ellipsoidal Bound on the Reachable Set} \label{sec:tightbound_reachable}
The expression of the reachable set in \eqref{eq:shapematrixstate2} indicates that the reachable set of the system state is a function of time, $k$, and composed of contributions from ellipsoidally bounded inputs, which are characterized with shape matrices, $\Sigma_{ij}$. Using the framework provided by this paper, the reachable set is bounded by a minimum trace ellipsoid $\mathcal{R}_{k}\subseteq\mathcal{E}(Q_k^*)$, where $\mathcal{E}(Q_k^*)$ is given by Lemma \ref{mintrace},
\begin{equation}
    \begin{aligned}
   Q_{k,\kappa}^*&=\left( \sum_{i=1}^r\sum_{j=\kappa+1}^{k-1} \sqrt{\mathrm{tr}(\Sigma_{ij})}\right)\left(\sum_{i=1}^r\sum_{j=1}^{k-1} \frac{\Sigma_{ij}}{\sqrt{\mathrm{tr}(\Sigma_{ij})}}\right),
    \end{aligned}
\end{equation}
and we omit the second subscript, $Q^*_k$, when $\kappa=0$.

It is intuitive that as the time horizon grows ($k$ becomes large), that the stability of the system matrices $A_k$ is required for the reachable set and the ellipsoidal bound to be bounded.


\begin{lemma}
If the state matrix ($A_k$) is stable for all $k\in\mathbb{N}$ and inputs $u_{i,k}$ are bounded, then the minimum trace ellipsoidal bound $\mathcal{E}(Q_k^*)$, is bounded.
\end{lemma}
\textit{Proof:}
Since the system is time-varying, the reachable set can potentially change dramatically from $k$ to $k+1$. Thus we fix $k$ and show that the sequence of historical ellipsoids $\{Q^*_{k,k-\kappa}\}$ as $\kappa$ approaches $k$ is Cauchy. As the limit of this sequence, $Q^*_k$ is then bounded. To do this we define the scalar sequence $S_\kappa$
\begin{equation}
    S_\kappa:=\sqrt{\mathbf{tr}(Q_{k,k-\kappa}^*)}=\sum_{i=1}^r\sum_{j=k-\kappa+1}^{k-1} \sqrt{\mathrm{tr}(\Sigma_{ij})}.
\end{equation}
Based on equation \eqref{eq:state_seq}, for $\delta_\kappa= S_{\kappa+1}-S_{\kappa}$,
\begin{equation}
    \lim_{\kappa \to k} \delta_\kappa=\sum_{i=1}^r\lim _{\kappa \to k}\sqrt{\mathrm{tr}(\Sigma_{i,k-\kappa})}=0
\end{equation}
and according to the ratio test,
\begin{equation*}
    \lim_{\kappa \to k} \frac{\delta_{\kappa+1}}{\delta_\kappa}= \lim_{\kappa \to k} \frac{\sum_{i=1}^r\sqrt{\mathrm{tr}(\Sigma_{i,k-\kappa-1})}}{\sum_{i=1}^r\sqrt{\mathrm{tr}(\Sigma_{i,k-\kappa})}}=\max_{j\in \overline{1,k}} (\rho[A_j]) < 1,
\end{equation*}
where $\rho[.]$ denotes the maximum eigenvalue. This implies $S_\kappa$ is Cauchy and convergent. Therefore, $\mathbf{tr}(Q_\kappa^*)=S_\kappa^2$ is Cauchy and convergent, which concludes $Q_k^*$ is bounded. 
\hspace*{\fill}~$\blacksquare$


\begin{figure}
    \centering
    \includegraphics[width=0.5\linewidth]{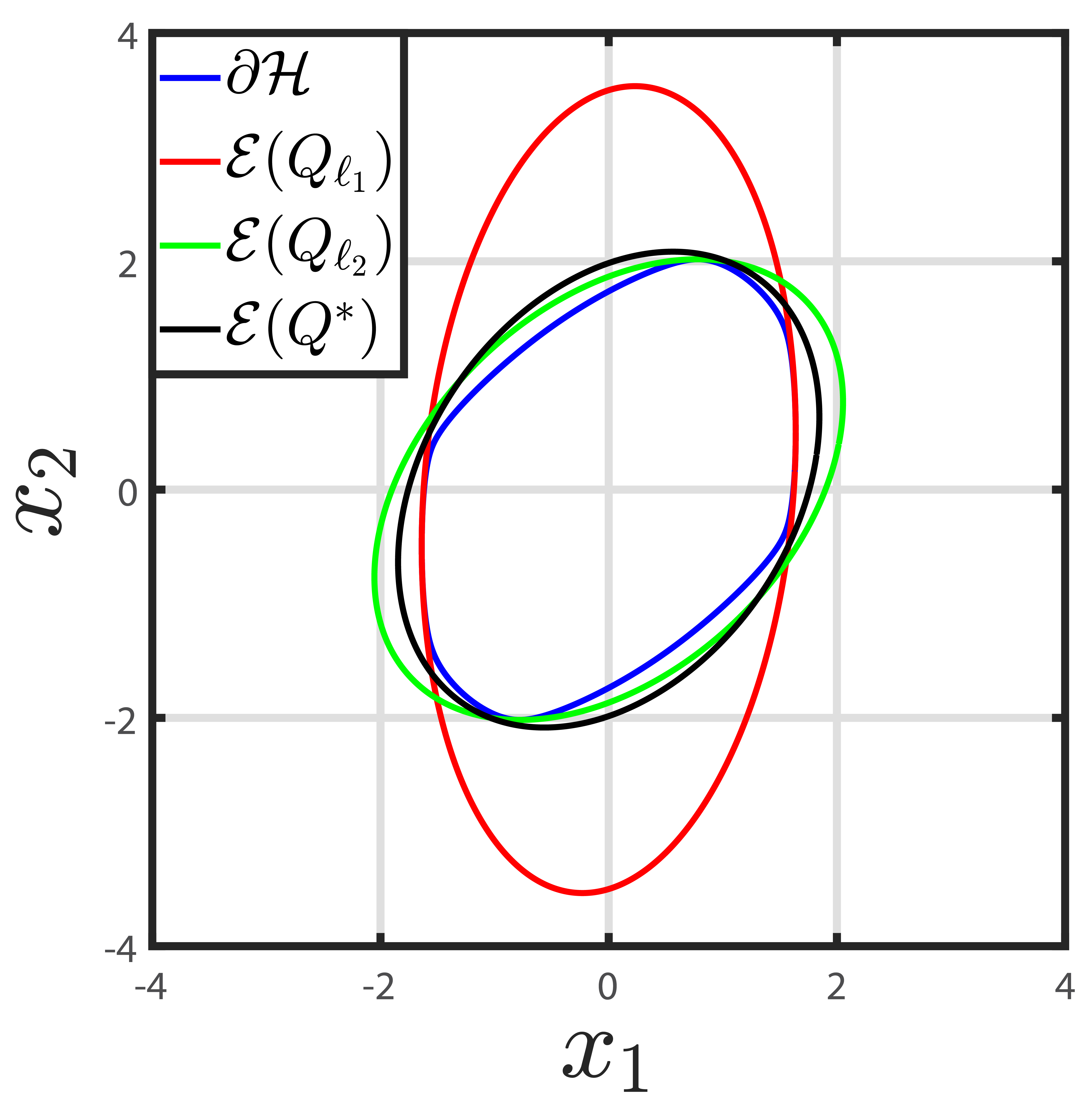}
    \caption{Exact boundary of geometric sum and a comparison between different ellipsoidal bounds.}
    \label{fig:comparison}
\end{figure}
\begin{figure}
    \centering
    \includegraphics[width=0.5\linewidth]{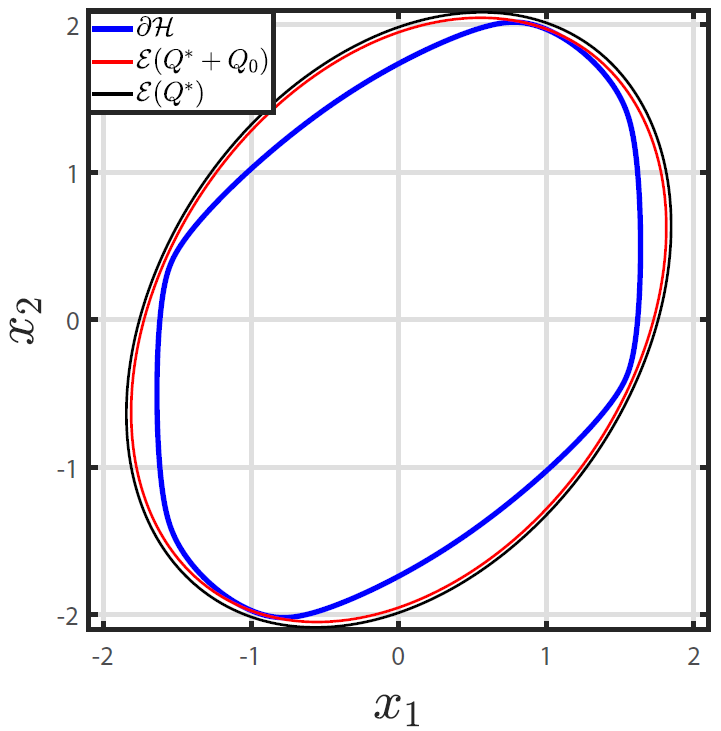}
    \caption{Selection of a nonzero $Q_0$ reduces the trace and produces a tangent ellipsoidal bound.}
    \label{fig:Q0}
\end{figure}
\begin{figure*}[t]
\centering
    \includegraphics[width=0.9\textwidth]{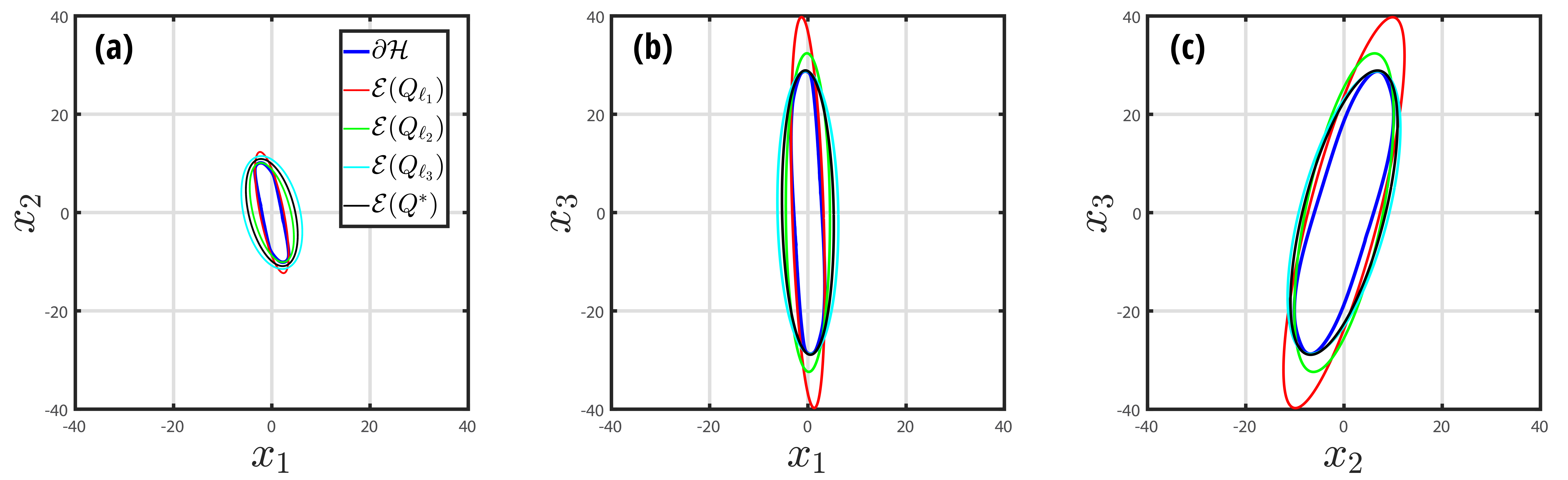}
  \caption{Exact boundary of the reachable set of a dynamical system and a comparison between different ellipsoidal bounds.}
  \label{fig:Reachsetprojection}
\end{figure*}

\section{Numerical Examples}
In this section we employ the developed tools to compute the minimum trace ellipsoidal bound over both a geometric sum and the reachable set of a dynamical system.

Consider the geometric sum of 4 ellipsoids with dimension 2 and shape matrices, 
{\small
$$
\begin{aligned}
Q_1&=\begin{bmatrix*}[r] 0.41&0.33\\0.33&0.31
\end{bmatrix*}, &Q_2=\begin{bmatrix*}[r]0.23&0.11\\0.11&0.06\end{bmatrix*},\\
Q_3&=\begin{bmatrix*}[r]0.17&-0.1\\-0.1&0.15\end{bmatrix*},\  &Q_4=\begin{bmatrix*}[r]0.01&0\\0&0.65\end{bmatrix*}.
\end{aligned}
$$}

We apply Theorem \ref{exactbound} to compute for the exact boundary of the geometric sum, plotted in blue in Fig. \ref{fig:comparison}. We then compute two tight outer ellipsoidal bounds using Corollary \ref{othershapematrix} corresponding to the standard unit vectors $\ell_1=\begin{bmatrix} 1& 0 \end{bmatrix}^T$ and $\ell_2=\begin{bmatrix} 0& 1 \end{bmatrix}^T$, producing $\mathcal{E}(Q_{\ell_1})$ and $\mathcal{E}(Q_{\ell_2})$ (red and green, respectively in Fig. \ref{fig:comparison}),
$$
Q_{\ell_1}=\begin{bmatrix*}[r] 2.68&0.81\\0.81&12.49\end{bmatrix*}, \quad Q_{\ell_2}=\begin{bmatrix*}[r] 4.22&1.57\\1.57&4.07\end{bmatrix*}.
$$
These ellipsoids are tangent to the exact boundary of the geometric sum. Lemma \ref{mintrace} returns the (black) minimum trace ellipsoid, $\mathcal{E}(Q^*)$,
$$
Q^*=\begin{bmatrix*}[r] 3.41 & 1.17\\1.17 & 4.34\end{bmatrix*}.
$$
As computed, $\mathcal{E}(Q^*)$ is not tangent to the geometric sum, however, it is visually the most representative ellipsoidal bound. We can numerically calculate a symmetric matrix $Q_0$ that reduces the value of the trace such that the resulting ellipsoidal bound becomes tangent to the exact geometric sum. For example,
$$
Q_0=\begin{bmatrix} -0.1193  & -0.0412\\ -0.0412 & -0.1521 \end{bmatrix}
$$
satisfies the conditions of Theorem \ref{subset}. The trace of  $\mathcal{E}(Q^*+Q_0)$ is smaller than the trace of $\mathcal{E}(Q^*)$ and is tangent to the geometric sum, see Fig. \ref{fig:Q0}.

In the second example, consider the LTI system
$$
x_{k+1}=F x_k + \nu_k, \quad F=\begin{bmatrix*}[r]0.67&0.35&-0.12\\-0.66&-0.55&0.41\\2.12&1.83&0.47\end{bmatrix*},
$$
where the noise $\nu_k \in \mathcal{E}(I_3)$, is spherically bounded.   

Here we compute the exact boundary of the reachable set, the minimum trace ellipsoidal bound, and, as before, compare them with 3 different tangent ellipsoids from Corollary \ref{othershapematrix}, $\mathcal{E}(Q_{\ell_1}),\mathcal{E}(Q_{\ell_2})$ and $\mathcal{E}(Q_{\ell_3})$ constructed by the standard unit vectors.
The settling time for this system is $k^*=120$, therefore, the reachable set can be closely approximated by the geometric sum of $120$ ellipsoids $\mathcal{E}(F^iF^{iT})$, $i\in\overline{1,120}$. The projections of the reachable set exact boundary and ellipsoidal bounds are shown in Fig.
\ref{fig:Reachsetprojection}.


\section{Conclusion}
In this paper we have generalized several results on characterizing the outer ellipsoidal bound on the geometric sum of a finite number of ellipsoids. Specifically, we introduce a symmetric matrix, $Q_0$, that can modify the existing formulas for these ellipsoidal outer bounds. The presence of this term exposes that the current formula for the minimum trace ellipsoid may not actually touch the reachable set and, therefore, cannot be the true minimum trace. Future work along these lines will hopefully reveal the relationship between the choice of this $Q_0$ term and the input ellipsoids ($Q_i$) and the choice of the coefficients ($p_{ij}$) that combine the shape matrices of the input ellipsoids to produce the outer ellipsoidal bounds.

\bibliographystyle{IEEEtran}
\bibliography{security}

\end{document}